\documentclass[12pt]{article}
\usepackage{amsfonts,amsmath,amssymb,amscd,amsthm}
\def\diam{\mathop{\fam0 diam}}

\long\def\comment#1\endcomment{}

\theoremstyle{theorem}
         \newtheorem{theorem}{Theorem}
         \newtheorem{lemma}[theorem]{Lemma}
         
         \newtheorem{proposition}[theorem]{Proposition}
         \newtheorem{conjecture}[theorem]{Conjecture}

\begin{document}

\centerline{\bf\uppercase{A two-page disproof of the Borsuk partition conjecture}}

\bigskip
\centerline{\bf A. Skopenkov
\footnote{Supported by a grant of Simons Foundation.
www.mccme.ru/\~{ }skopenko}
}

\bigskip
\small
{\bf Abstract.} It is presented the simplest known disproof of the Borsuk conjecture stating that
{\it if a bounded subset of $n$-dimensional Euclidean space contains more than $n$ points,
then the subset can be partitioned into $n+1$ nonempty parts of smaller diameter.}
The argument is due to N. Alon and is a remarkable application of combinatorics and algebra to geometry.
This note is purely expository and is accessible for
students.

\normalsize

\bigskip
\begin{theorem}[Borsuk]\label{bor}
If a bounded subset of the plane contains more than two points, then the subset
 can be partitioned into three nonempty parts of smaller diameter.
\footnote {{\it Hint to the proof.}
First, using the `considerations of continuity ', prove that any plane figure of diameter 1 can be
included in a regular hexagon whose inscribed circle has diameter 1.
Then prove that, although the diameter of the obtained regular hexagon is greater than 1, the hexagon
can be cut into three parts of diameter less than 1.
Cf. [Y10]. }
\end{theorem}

The {\it diameter} of a nonempty subset of the plane is the greatest distance between points of the subset.
(More precisely, the supremum of such distances).
 A subset of the plane is called {\it bounded}, if its diameter is finite.

Borsuk proposed the following generalization,
which for years was one of the most intriguing conjectures in combinatorial geometry.

A {\it point}
$x=(x_1,\dots,x_n)$ of $n$-dimensional Euclidean space is an ordered set of $n$ numbers.
The {\it distance} between points $x=(x_1,\dots,x_n)$ and $y=(y_1,\dots,y_n)$ is defined by
$$|x,y|:=\sqrt{(x_1-y_1)^2+\dots+(x_n-y_n)^2}.$$
The {\it diameter} and the property of being {\it bounded} for a subset of $n$-dimensional
Euclidean space is defined the same way as for a subset of the plane.

\begin{conjecture}[Borsuk]\label{borcon}
If a bounded subset of $n$-dimensional Euclidean space contains more than $n$ points,
then the subset can be partitioned into $n+1$ nonempty parts of smaller diameter.
\end{conjecture}

In 1993 J. Kan and G. Kalai, following the ideas Boltyanskii, Erd\"os and Larman on the application of
combinatorics to construct a counterexample, found a counterexample to the Borsuk conjecture [KK93, S96].
Background is described in detail in [AZ04, R06].

 \begin{theorem}\label{kk}
There exist $n$ and a bounded subset of $n$-dimensional Euclidean space such that the subset
contains more than $n$ points and can not be partitioned into $n+1$ parts of smaller diameter.
  \end{theorem}

We present the simplest known proof due to N. Alon, cf. [N94, G99, S99, AZ04, R04, R06].
(Note that other proofs give stronger results.)
This is a wonderful example of an important result in modern mathematics whose proof does not require a full understanding of a special six-month university course (after a two-year ordinary course).
Main steps of the proof are Estimation (Lemma \ref{estlem}) and Linear Independence (Lemma \ref{lin}).


By $|X|$ we denote the number of elements in a set $X$.
The {\it scalar product} of points (=vectors) $x=(x_1,\dots,x_n)$ and $y=(y_1,\dots,y_n)$ is defined by
$x\cdot y:=x_1y_1+\dots+x_ny_n$.
Vectors $x$ and $y$ are called {\it orthogonal} if $x\cdot y=0$.

\smallskip
 {\it Proof of Theorem \ref{kk}.}
Denote
$$ M = \{(x_1, \dots, x_n) | \ x_1 = 1, \ x_k \in \{1, -1 \} \mbox{ and }x_2\cdot\ldots\cdot x_n=1\}.$$
 The vertex of $n^2$-dimensional cube is the ordered set of $n^2$ plus or minus units.
It can be thought of as a table $n\times n$.
 We assign to each point $x=(x_1, \dots, x_n) \in M$ table $fx$ defined by the formula $(fx)_{ij}: = x_ix_j$.
For example
$$f(1, -1, -1)=\left (\begin{matrix} 1 & -1 & -1 \\ -1 & 1 & 1 \\ -1 & 1 & 1 \end{matrix} \right).$$
Let us prove that a counterexample to Borsuk's conjecture is the $f$-image of $M$ for a sufficiently large
prime number $p$ and $n = 4p$.

Take $x,y \in M$. Then $(x_ix_j-y_iy_j)^2=(1-x_iy_ix_jy_j)^2$.
Denote by $a$ the number of subscripts $i$ for which $x_i=y_i$.
Then $x_iy_i=1$ for $a$ subscripts $i$ and $x_iy_i =-1$ for $n-a$ subscripts $i$.
Therefore $|fx,fy|^2=4a(n-a)$.
This is maximal for $a=n/2$.
Hence the condition $ | fx, fy | = \diam fM $ is equivalent to $a=n/2$ and is equivalent to $x\cdot y=0$.

Thus if the set $fM$ is partitioned into $k$ parts $Z_1,\dots,Z_k$ of smaller diameter, then in each
$f^{-1}Z_i$ no two vectors are orthogonal.
Since $x_1=1$ for each $x\in M$, the map $f$ is injective.
Hence  $|Z_i|=|f^{-1}Z_i|$.
Now the theorem is implied by the following Estimation (Lemma \ref{estlem}) and Proposition \ref{estim}.
QED

 \begin{lemma}[Estimation]\label{estlem}
Assume that $p$ is a prime (not necessarily large)
$n=4p$, $A\subset M$ and no two vectors of $A$ are orthogonal. Then
$$|A|\le\alpha(n): ={n-1\choose0} + {n-1\choose1} + \dots + {n-1\choose\frac n4-1}.$$
  \end{lemma}

\begin{proposition} \label{estim}
$\alpha(n)(n^2+1)<|M|=2^{n-2}$ for sufficiently large $n$.
\end{proposition}

{\it Proof.}
For sufficiently large $n$ and each $s,k\in\{1,2,\dots,n/4\}$:
$$\frac{5n}4>5k-s-2 \quad \Rightarrow \quad \frac{nk-s +1}{\frac n4 + ks}>\frac32 \quad \Rightarrow $$
$$ \Rightarrow \quad \frac{{n-1\choose\frac n4 + k-1}}{{n-1\choose k-1}} =
\frac{(nk)(nk-1)\dots(nk-\frac n4 +1)}{(\frac n4+k-1)(\frac n4+k-2)\dots k}>(\frac32)^{n/4}>n^2 \quad \Rightarrow $$
$$ \Rightarrow \quad \alpha(n)(n^2+1)<{n-1\choose0}+{n-1\choose1}+\dots+{n-1\choose \frac n2-1}=2^{n-2}.\qquad QED$$

It remains to prove the Estimation (Lemma \ref{estlem}).
In its proof one can forget about the construction of $f$.
The following proposition is clear.

 \begin{proposition}\label{div}
For a prime integer $p$ and an integer $ t $ the number
$$G(t):=(t-1)(t-2)\dots(t-p +1)$$
is divisible by $p$ if and only if $t$ is not divisible by $p$.
\end{proposition}

{\it Rational linear combinations} of polynomials $F_1,\dots,F_s$ is a polynomial $\lambda_1F_1+\dots+\lambda_sF_s$ with rational $\lambda_1,\dots,\lambda_s$.
For example, the polynomial $x_2$ is a rational linear combination of polynomials $2x_1$, 1 and $x_1+x_2$.

Polynomials are called {\it linearly independent} if each their rational linear combination with $\lambda_k$'s not all  zeros, is non-zero.
For example, $n$ polynomials $1,x_2,x_3,\dots,x_n$ are linearly independent.

A polynomial with rational coefficients in $n-1$ variables $x_2,\dots,x_n$ is called {\it square-free and of degree less than $n/4$}
if it is a rational linear combination of polynomials
$$(*)\qquad x_ {i_1}\cdot\ldots\cdot x_ {i_s}, \quad \text{where} \quad s\in\{0,\dots,n/4-1\} \quad \text{and}
\quad i_1,\dots,i_s\in\{2,3,\dots,n\}\text{ are distinct}.$$
The Estimation (Lemma \ref{estlem}) is implied by the Linear Independence (Lemma \ref{lin})
and Proposition \ref{mains} below.

\begin{lemma}[Linear Independence]\label{lin}
  Assume that $p$ is a prime integer, $n=4p$, \ $A\subset M$ and no two vectors of $A$ are orthogonal.
Take a vector $a\in A$.
Open brackets in the product $G(a \cdot(1,x_2,\dots,x_n))$.
In each of the resulting monomials for each $i$ replace $x_i^2$ by 1 as long as possible.
The resulting polynomial is denoted by $F_a(x_2,\dots,x_n)$.
Then each polynomial $F_a(x_2,\dots,x_n)$, $a\in A$, is square-free and of degree less than $n/4$; these polynomials are linearly independent.
  \end{lemma}

\begin{proposition} \label{mains}
A linearly independent family of square-free polynomials in $x_2,\dots,x_n$ of degree less than $n/4$ contains at most $\alpha(n)$ polynomials.
  \end{proposition}

 {\it Proof of the Linear Independence Lemma.} The `square-free' and the `degree' properties are obvious.
Let us prove the linear independence. Suppose to the contrary that
$\lambda_1F_{a_1}+\dots+\lambda_sF_{a_s}=0$
for some $a_1,\dots,a_s\in A$ and rational $\lambda_1,\dots,\lambda_s$, not all $\lambda_k$ being zero.
Here $a_1,\dots,a_s$ are vectors not coordinates.
We may assume that $\lambda_1,\dots,\lambda_s$ are integers, not all of them divisible by $p$. 
(Otherwise multiply this equation by the product of denominators and then divide by the greatest common divisor.)
Without loss of generality $\lambda_1$ is not divisible by $p$.
In the obtained equality take $x_2=(a_1)_2,\dots,x_n=(a_1)_n$.

Since $a_1\cdot a_1=n=4p$, by Proposition \ref{div} $\lambda_1F_{a_1}$ is not divisible by $p$.

Since $n$ is divisible by 4 and for each $x,y\in A$ the numbers of minus units in $x$ and in $y$ are even,
$x\cdot y$ is divisible by 4.
Therefore $x\cdot y\not\in\{\pm p,\pm2p,\pm3p\} $.
Since $x\cdot y\ne0$, the product $x\cdot y$ is not divisible by $p$.
Hence by Proposition \ref{div}  $\lambda_kF_{a_k}$ is divisible by $p$ for each $k>1$.
A contradiction.
QED

\smallskip
{\it Sketch of a proof of Proposition \ref{mains}.}
 Denote by $Q_1,\dots,Q_{\alpha(n)}$ the polynomials (*) and
by $F_1,\dots,F_k$ given linearly independent family.
Take the $k\times\alpha(n)$ table of rational numbers $\lambda_{ij}$ for which $F_i=\sum_j\lambda_{ij}Q_j$ when $i=1,\dots,k$.
The family of polynomials obtained from the family $F_1,\dots,F_k $ substituting $F_i$ by $F_i+\lambda F_j$, $j\ne i$,
is linearly independent.
Using such substitutions  and permutations of polynomials $Q_1,\dots,Q_{\alpha (n)}$ we can change the
$k\times\alpha(n)$ table to the `upper triangular' form.
Since the new table has no zero lines,  $k\le\alpha(n)$.
QED

\bigskip
{\bf Acknowledgements.}

I would like to acknowledge N.P. Dolbilin and A.M. Raigorodskii from whom I learned
counterexamples to the Borsuk Conjecture, students of Kolmogorov College and of high-school N57 in Moscow
who have learned from these counterexamples from me, and M.B. Akhmedov, V.N. Dubrovsky and A.D. Rukhovich for useful discussions.

\bigskip
{\bf References.}

[AZ04] M. Aigner, G. Ziegler, Proofs from the Book, Springer, 2004.

[G99] M.L. Gerver, On partition of sets into parts of smaller diameter (in Russian),
Mat. Prosveschenie, 3 (1999) 168-183.

[KK93] J. Kahn and G. Kalai, A counterexample to Borsuk's conjecture, Bull. AMS, 29:1 (1993) 60--62.

[N94] A. Nilli, On Borsuk's problem, Contemp. Math., 178 (1994) 209--210.

[R04] A. M. Raigorodskii, The Borsuk partition problem: the seventieth anniversary,
Math. Intelligencer, 26:3 (2004) 4--12.

[R06] A. M. Raigorodskii, The Borsuk problem (in Russian), MCCME, Moscow, 2006.

[S96] A. Skopenkov, The Borsuk problem, Quantum, 7:1 (1996) 16--21, 63.

[S99] A. Skopenkov, $N$-dimensional cube, polynomials and solution of the Borsuk
problem (in Russian), Mat. Prosveschenie, 3 (1999) 184--188.
http://arxiv.org/abs/0712.4009, v1.

[Y10] Dian Yang, An elementary proof of Borsuk theorem,
http://arxiv.org/abs/1010.1990.

\end{document}


\centerline{\bf\uppercase{Короткое опровержение гипотезы Борсука}
}

\bigskip
\centerline{\bf А. Скопенков
\footnote{Поддержан грантом фонда Саймонса. www.mccme.ru/~skopenkov}
}

\bigskip
\small
{\bf Аннотация.} Приводится простейшее из известных опровержений следующей гипотезы Борсука:
{\it любое ограниченное подмножество $n$-мерного евклидова пространства, содержащее более $n$ точек,
можно разбить на $n+1$ непустых частей меньшего диаметра. }
Доказательство принадлежит Н. Алону и является замечательным приложением комбинаторики и алгебры к геометрии.
Эта методическая заметка доступна студентам, старшеклассникам и учителям, интересующимся математикой.

\normalsize
\smallskip
\begin{theorem}[Борсук]\label{bor}
Любое ограниченное подмножество плоскости, в котором более двух точек,
можно разбить на три непустые части меньшего диаметра.
\footnote{{\it Указание к доказательству.}
Сначала, используя `соображения непрерывности', докажите, что любую плоскую
фигуру диаметра 1 можно заключить в правильный шестиугольник, диаметр
вписанной окружности которого равен 1. Затем докажите, что хотя диаметр
полученного правильного шестиугольника больше 1, его можно разрезать на
три части диаметра меньше 1.  Ср. [Y10].}
\end{theorem}

 {\it Диаметром} непустого подмножества плоскости называется наибольшее расстояние между его точками
(точнее, супремум таких расстояний).
Подмножество плоскости называется {\it ограниченным}, если его диаметр конечен.

Борсук предложил следующее обобщение своего результата,
которое долгие годы было одной из наиболее интригующих проблем комбинаторной геометрии.

{\it Точкой}
$x=(x_1,\dots,x_n)$
$n$-мерного евклидова пространства называется упорядоченный набор $n$ чисел.
{\it Расстояние} между точками $x=(x_1,\dots,x_n)$ и $y=(y_1,\dots,y_n)$ определяется формулой
$$|x,y|:=\sqrt{(x_1-y_1)^2+\dots+(x_n-y_n)^2}.$$
{\it Диаметр} и {\it ограниченность} подмножества $n$-мерного евклидова пространства определяется точно так же,
как и в случае плоскости.

 \begin{conjecture}[Борсук]\label{borcon}
Любое ограниченное подмножество $n$-мерного евклидова пространства, содержащее более $n$ точек, можно разбить на $n+1$ непустых частей меньшего диаметра.
\end{conjecture}

В 1993 Д. Кан и Дж. Калаи, следуя идеям Болтянского, Эрдеша и Лармана о применении комбинаторики для построения
контрпримера, нашли контрпример к гипотезе Борсука [KK93, S96].
Подробно история вопроса описана в [AZ04, R06].




\begin{theorem}\label{kk}
Существует $n$ и ограниченное подмножество $n$-мерного евклидова пространства,
содержащее более $n$ точек и которое невозможно разбить на $n+1$ часть меньшего диаметра.
\end{theorem}

 Мы приведем простейшее из известных доказательств, принадлежащее Н. Алону, ср.
[N94, G99, S99, AZ04, R04, R06].
(При этом другие доказательства дают более сильные результаты.)
Это удивительный пример важного результата в современной математике, не требующего для полного понимания
полугодового специального университетского курса (после двухгодового обязательного курса).
Более простые применения аналогичных алгебраических соображений в комбинаторике можно найти в [IKS, R].

Через $|X|$ обозначается число элементов в множестве $X$.
{\it Скалярное произведение} векторов $x=(x_1,\dots,x_n)$ и $y=(y_1,\dots,y_n)$ определяется
как $x\cdot y:=x_1y_1+\dots+x_ny_n$.
Векторы $x$ и $y$ называются {\it ортогональными}, если $x\cdot y=0$.

\smallskip
{\it Доказательство теоремы \ref{kk}.}
Обозначим
$$M=\{(x_1,\dots,x_n)|\ x_1=1,\ x_k\in\{1,-1\}\mbox{
 и среди }x_2,\dots,x_n\mbox{ число минус единиц четно}\}.$$
Вершина $n^2$-мерного куба --- набор длины $n^2$ из плюс или минус единиц.
Его удобно представлять себе как таблицу $n\times n$.
Поставим в соответствие каждой точке $x=(x_1,\dots,x_n)\in M$
таблицу $fx$, определенную формулой $(fx)_{ij}:=x_ix_j$.
Например,
$$f(1,-1,-1)=\left(\begin{matrix} 1 & -1 & -1 \\ -1 & 1 & 1 \\ -1 & 1 & 1 \end{matrix}\right)$$
Докажем, что контрпримером к гипотезе Борсука является $f$-образ множества $M$ для достаточно большого простого числа $p$ и $n=4p$.



Пусть $x,y\in M$. Тогда $(x_ix_j-y_iy_j)^2=(1-x_iy_ix_jy_j)^2$.
Обозначим через $a$ количество индексов $i$, для которых $x_i=y_i$.
Тогда $x_iy_i=1$ для $a$ индексов $i$  и $x_iy_i=-1$ для $n-a$ индексов $i$.
Поэтому $|fx,fy|^2=4a(n-a)$.
Это выражение максимально при $a=n/2$.
Значит, условие $|fx,fy|=\diam fM$ равносильно условию $a=n/2$ и равносильно условию $x\cdot y=0$.

Поэтому если множество $fM$ разбито на $k$ частей $Z_1,\dots,Z_k$ меньшего диаметра, то
в каждом $f^{-1}Z_i$ никакие два вектора не ортогональны.
Так как $x_1=1$ для любого $x\in M$, то $f$ инъективно.
Значит, $|Z_i|=|f^{-1}Z_i|$.
Теперь
теорема вытекает из следующих леммы \ref{estlem} об оценке и утверждения \ref{estim}.
QED

 \begin{lemma}[Оценка]\label{estlem}
 Пусть $p$ -- простое (не обязательно большое),
$n=4p$, $A\subset M$ и никакие два вектора из $A$ не ортогональны. Тогда
$$|A|\le\alpha(n):=C_{n-1}^0+C_{n-1}^1+\dots+C_{n-1}^{\frac n4-1}.$$
\end{lemma}

\begin{proposition} \label{estim}
 $\alpha(n)(n^2+1)<|M|=2^{n-2}$ для достаточно больших $n$.
 \end{proposition}

 {\it Доказательство.}
Для достаточно больших $n$ и любых $s,k\in\{1,2,\dots,n/4\}$:
$$\frac{5n}4>5k-s-2\quad\Rightarrow\quad \frac{n-k-s+1}{\frac n4+k-s}>\frac32\quad\Rightarrow$$
$$\Rightarrow\quad\frac{C_{n-1}^{\frac n4+k-1}}{C_{n-1}^{k-1}}=
\frac{(n-k)(n-k-1)\dots(n-k-\frac n4+1)}{(\frac n4 +k-1)(\frac n4 +k-2)\dots k}>(\frac32)^{n/4}>n^2\quad\Rightarrow$$
$$\Rightarrow\quad \alpha(n)(n^2+1)<C_{n-1}^0+C_{n-1}^1+\dots+C_{n-1}^{\frac n2-1}=2^{n-2}.\qquad QED$$

 Осталось доказать лемму об оценке.
При ее доказательстве можно забыть про конструкцию отображения $f$.
Следующее утверждение очевидно.

 \begin{proposition}\label{div}
 Для простого $p$ и целого $t$ число
$$G(t):=(t-1)(t-2)\dots(t-p+1)$$
делится на $p$ тогда и только тогда, когда $t$ не делится на $p$.
\end{proposition}

{\it Рациональной линейной комбинацией} многочленов $F_1,\dots,F_s$ называется
любой многочлен $\lambda_1F_1+\dots+\lambda_sF_s$ с рациональными $\lambda_1,\dots,\lambda_s$.
Например, многочлен $x_2$ является рациональной линейной комбинацией многочленов $2x_1$, $1$ и $x_1+x_2$.

Многочлены называются {\it линейно независимыми}, если любая их рациональная линейная комбинация, в которой
не все $\lambda_k$ нулевые, не равна нулю.
Например, $n$ многочленов $1,x_2,x_3,\dots,x_n$ являются линейно независимыми.

Многочлен с рациональными коэффициентами от $n-1$ переменной $x_2,\dots,x_n$ называется {\it степени менее $n/4$ и свободным от квадратов},
если он является рациональной линейной комбинацией многочленов
$$(*)\quad x_{i_1}\cdot\ldots\cdot x_{i_s},\quad\text{где}\quad s=0,\dots,p-1\quad\text{и}\quad i_1,\dots,i_s\quad
\text{ --- различные числа от $2$ до $n$}.$$
Лемма об оценке вытекает из нижеследующих леммы \ref{lin} о линейной независимости и утверждения \ref{mains}.

\begin{lemma}[Линейная независимость]\label{lin}
 Пусть $p$ простое, $n=4p$, \ $A\subset M$ и никакие два вектора из $A$ не ортогональны.
Возьмем вектор $a\in A$.
Раскроем скобки в произведении $G(a\cdot (1,x_2,\dots,x_n))$.
В каждом из полученных одночленов для каждого $i$ будем заменять $x_i^2$ на $1$, пока это возможно.
Полученный многочлен обозначим $F_a(x_2,\dots,x_n)$.
Тогда каждый многочлен $F_a(x_2,\dots,x_n)$, $a\in A$, степени меньше $n/4$ и свободен от квадратов; эти многочлены линейно независимы.
 \end{lemma}

 \begin{proposition} \label{mains}
Любое линейно независимое семейство многочленов от $x_2,\dots,x_n$ степени менее $n/4$ и свободных от квадратов, содержит не более $\alpha(n)$ многочленов.
\end{proposition}

 {\it Доказательство леммы о линейной независимости.}
 Утверждения о степени и о свободе от квадратов очевидны.
Докажем линейную независимость.
Пусть, напротив,
$\lambda_1F_{a_1}+\dots+\lambda_s F_{a_s}=0$
для некоторых
$a_1,\dots,a_s\in A$
и рациональных
$\lambda_1,\dots,\lambda_s$, причем не все $\lambda_k$ нулевые.
Здесь $a_1,\dots,a_s$ --- векторы, а не координаты.
Можно считать, что $\lambda_1,\dots,\lambda_s$ целые (иначе умножим это равенство на произведение их знаменателей).
Можно считать, что не все они делятся на $p$ (иначе поделим это равенство на их наибольший общий делитель).
Не уменьшая общности, $\lambda_1$ не делится на $p$.
Подставим в полученное равенство значения $x_2=(a_1)_2,\dots,x_n=(a_1)_n$.

Из $a_1\cdot a_1=n=4p$ и утверждения \ref{div} вытекает, что $\lambda_1F_{a_1}$ не делится на $p$.

Так как $n$ делится на 4 и для любых $x,y\in M$ число минус единиц в $x$ и в $y$ нечетно, $x\cdot y$ делится на 4.
Поэтому $x\cdot y\not\in\{\pm p,\pm2p,\pm3p\}$.
Так как $x\cdot y\ne0$, то $x\cdot y$ не делится на $p$.
Значит, по утверждению \ref{div} $\lambda_kF_{a_k}$ делится на $p$ при любом $k>1$.
Противоречие.
QED

\smallskip
{\it Набросок доказательства утверждения \ref{mains}.}
Обозначим через $Q_1,\dots,Q_{\alpha(n)}$ семейство многочленов (*) и
через $F_1,\dots,F_k$ данное линейно независимое семейство.
Возьмем таблицу $k\times\alpha(n)$ рациональных чисел $\lambda_{ij}$,
для которых $F_i=\sum_j\lambda_{ij}Q_j$ при любом $i=1,\dots,k$.
Семейство многочленов, полученное из семейства $F_1,\dots,F_k$ заменой $F_i$ на $F_i+\lambda F_j$, $j\ne i$,
линейно независимо.
Такими заменами и перестановками многочленов $Q_1,\dots,Q_{\alpha(n)}$ можно провести
рассматриваемую таблицу $k\times\alpha(n)$ к `верхнетреугольному' виду.
Так как в новой таблице нет нулевой строки, то $k\le\alpha(n)$.
QED



\bigskip
{\bf Благодарности.}
Благодарю Н.П. Долбилина и А.М. Райгородского, от которых я узнал
контрпримеры к гипотезе Борсука, учеников физ.-мат. школы им.
А.Н. Колмогорова и школы N57 г. Москвы, которые узнали эти контрпримеры от
меня, а также М.Б. Ахмедова, В.Н. Дубровского и А.Д. Руховича за полезные обсуждения.

\bigskip
{\bf Литература.}

[AZ04] M. Aigner, G. Ziegler, Proofs from the Book, Springer, 2004.


[G99] М. Л. Гервер, О разбиении множеств на части меньшего диаметра: теоремы и контрпримеры,
Мат. Просвещение, 3 (1999), 168-183.

[IKRS] М. Ильинский, А. Купавский, А. Райгородский и А. Скопенков,
Дискретный анализ для математиков и программистов (подборка задач),
представлено к публикации.

[KK93] J. Kahn and G. Kalai, A counterexample to Borsuk's conjecture, Bull. AMS, 29:1 (1993) 60--62.

[N94] A. Nilli, On Borsuk's problem, Contemp. Math., 178 (1994) 209--210.

[R04] A. M. Raigorodskii, The Borsuk partition problem: the seventieth anniversary,
Math. Intelligencer, 26:3 (2004) 4--12.

[R06] А. М. Райгородский, Проблема Борсука, МЦНМО, Москва, 2006.

[R] А.М. Райгородский, Линейно-алгебраический метод в комбинаторике, Москва, МЦНМО.

[S96] A. Skopenkov, The Borsuk problem, Quantum, 7:1 (1996) 16--21, 63.

[S99] А. Скопенков, $N$-мерный куб, многочлены и решение проблемы Борсука,
Мат. Просвещение, 3 (1999) 184--188.
http://arxiv.org/abs/0712.4009, v1

[Y10] Dian Yang, An elementary proof of Borsuk theorem,
http://arxiv.org/abs/1010.1990.